\documentclass[a4paper]{article}

\usepackage[english]{babel}
\usepackage[utf8x]{inputenc}
\usepackage[T1]{fontenc}
\usepackage{setspace}
\usepackage{multicol}
\usepackage{graphicx}
\usepackage{float}
\graphicspath{ {images/} }

\usepackage[a4paper,top=3cm,bottom=2cm,left=3cm,right=3cm,marginparwidth=1.75cm]{geometry}

\usepackage[affil-it]{authblk}
\usepackage{amsmath}
\usepackage{graphicx}
\usepackage[colorinlistoftodos]{todonotes}
\usepackage[colorlinks=true, allcolors=blue]{hyperref}
\begin{document}

\title{\Huge \textbf{Minimal Representations of Natural Numbers Under a Set of Operators}}
\author{\LARGE \textbf{Akshunna Shaurya Dogra}%
 \thanks{Email: \texttt{adogra@nyu.edu}}}
\affil{Courant Institute of Mathematical Sciences,\\ New York University}

\maketitle

\begin{abstract}\small
This paper studies the minimal length representation of the natural
numbers. Let $O$ be a fixed set of integer-valued functions (primarily hyperoperations). For each $n$, what
is the shortest way of expressing $n$ as a combinations of functions in $O$ to
the constant $1$? For example, if $O$ contains the two functions $S'x'$ (successor
of $'x'$) and $*'x''y'$ (x times y) then the shortest representation of 12 is $'*SSS1SS1'$,
with 8 symbols. This is taken to mean that 8 is complexity of 12 under $O$.

We make a study of such minimal representations and complexities in this paper, proving and/or rightly predicting bounds on complexities, discussing some relevant patterns in the complexities and minimal representations of the natural numbers and listing the results gleaned from computational analysis. Computationally, the first 4.5 million natural numbers were probed to verify our mathematically obtained results. Due to the finiteness of the problem, we used the method of exhaustion of possibilities to state some other results as well.
\end{abstract}
\pagebreak
\section{Introduction and some general results}
Let $O$ be a set of functions and constant symbols over the natural numbers. Assume that each function has a fixed arity, so that expressions can be written unambiguously in prefix notation. \\
\\
  \textbf{Definition 1.1: A \textit{term} is an expression over $O$ in prefix notation; we will write these as strings delimited by single quotes. We will denote the value of term} $'a'$ \textbf{as} \textit{v(a)} \\ 
\linebreak
\textbf{Definition 1.2: A \textit{representation}} $'r'$ \textbf{of a number} \textit{n} \textbf{over $O$ is a term using the symbols in $O$ that evaluates to} \textit{n.} \\ 
\linebreak
\textbf{Definition 1.3: The \textit{length} of} $'r'$, \textbf{denoted as} \textit{\#r}\textbf{, is the number of symbols in} \textit{r}.\\
\linebreak
\textbf{Definition 1.4: A \textit{minimal representation}
of $n$ is a representation of $n$ with the fewest symbols. The complexity of $n$
over $O$, denoted $c_{O}(n)$ is the length of a minimal-length representation.} \textbf{Therefore,  }$c_{O}(n) = min\{\#r | v(r) = n\}$.
\\
\\
\textbf{Definition 1.5: $M_{O}(k)$ is the representation of length $k$ of maximal value. That is, $\forall$ $'r'$, $\#r = k\implies v(r) \leq v(M_{O}(k)$}\\ 
\textbf{Example 1.1:} Let $O = \{1, S, *\}$($S$ is the successor function).  Then $'*SS1SSS1'$, $'*S1*S1SS1'$, and $'SS*S1SSSS1'$, $'SSSSSSSSSSS1'$ are all representations of 12, with lengths 8, 9, 10 and 12, respectively. $v(*SS1SSS1)$ = 12. $'*SS1SSS1'$ is the minimal representation under $'O'$, thus the complexity being 8 (each operation is followed by its left and right operand. Therefore, $'+SS1*SS1 \wedge S1SS1'$ represents $3 + 3(2^3)$ and so on). As 12 is the largest number possible from 8 symbols under $O$, $'*SS1SSS1'$ is also $'M_O(k)'$ for $k=8$ 
\\
\\
\textbf{\textit{Theorem 1.1: Let $O = \{1, S\}$. Then,}} $\forall n$, \textit{c\textsubscript{O}}$(n) = n$.
\\
\textit{\textbf{Proof}:} The only representation of $n$ consists of $n − 1$ occurrences of $S$ applied to $1$; and this has length $n$.
\\
\textbf{Corollary: If $\{1, S\} \subset O$}, \textit{c\textsubscript{O}(n)} $\leq n$.\\
\linebreak
\textbf{Definition 1.6: A number $n$ is called \textit{irreducible} in $O$ if $\{1, S\}\subset O$ \textbf{\&}} \textit{c\textsubscript{O}(n) = n}.
\linebreak
\\
\textbf{\textit{Theorem 1.2: Let $O = \{1, S, +\}$. Then,}} $\forall n, c\textsubscript{O}(n) = n$\textbf{.}
\\
\textit{\textbf{Proof:}} (By induction on $\#r$.) For any representation $'r', \#r ≥ v(r)$.
\\
Clearly it holds if $\#r = 1$. Assume that the statement holds for all terms $u$ where $\#u < n$. Let $r$ be a representation of length $n$. 
\begin{enumerate}
\item $'r'$ has the form $S'x'$, where \textit{\#x = n−1}. By the induction hypothesis \textit{v(x)} $\leq$ \textit{n−1}, so \textit{v(r)} $\leq$ \textit{n}.
\item $'r'$ has the form $+'x''y'$, where \textit{\#x+\#y = n−1}. By the induction hypothesis, \textit{v(x)+v(y)} $\leq$ \textit{n−1} so \textit{v(r)} $\leq$ \textit{n}.
\end{enumerate}
Finally, if $'r'$ is the term consisting of $(n − 1)$ $S$ followed by $1$, then $v(r) = \#r = n$.\\
\\
\textbf{\textit{Theorem 1.3: Let $S \in O$. Let $n = v(M_{O}(k))$. Then $c_{O}(n) = k$}}
\\
\textbf{\textit{Proof:}} By construction, \textit{c\textsubscript{O}(n) $\leq$ k}. If \textit{c\textsubscript{O}(n) $= k' <$ k}, I can use \textit{($k-k'$) S} on $n'$ to make \textit{$n' > n$}, using \textit{k} symbols in total, giving us a contradiction. Therefore \textit{c\textsubscript{O}(n) $= k$.}
\\
\textbf{Corollary: Let $S \in O$.
For every length $k$, there exists a $n$ such that $k=c_{O}(n)$ (namely $n=v(M_{O}(k))$}\\
\\
\textbf{Definition 1.7: A number $n_u$ is \textit{ugly} under $O$, if $\forall m \neq n_u$,} \textit{c\textsubscript{O}(m)} $=$ \textit{c\textsubscript{O}($n_u$)} $ = k\implies m > n_u$\\
\\
\textbf{Definition 1.8: A number $n_e$ is \textit{efficient} in $O$, if $\forall m \neq n_e$,} \textit{c\textsubscript{O}(m)} $=$ \textit{c\textsubscript{O}($n_e$)} $=k \implies n_e > m$.
\\
\\
\textbf{\textit{Theorem 1.4: $\forall O, \exists$ infinitely many $n$ such that $c_{O}(n) \geq \log(n)/\log(|O|)$. More precisely, for $\epsilon > 0$, as $x \to\infty$, number of integers such that $c_{O}(n) < \log(n)/\log(|O|)$ is $\mathcal{O}(x^{1-\epsilon})$.}} \\
\textbf{\textit{Proof:}} Let $|O| = a, (1 - \epsilon)/\log(a) = z$. Let $l(k)$ be the number of possible terms evaluating to a natural number using \textit{k} symbols. Therefore, number of possible representations that evaluate to any integer using \textit{k} symbols from $O$ is $< a\textsuperscript{k}$.\\
Let $E(x) = \{n |n \leq x, c\textsubscript{O}(n) < z\log(n)\}$. \\
Say $k = max\{c\textsubscript{O}(n)| n \in E(x)\} \implies k < z\log(x)$.\\
Therefore, $|E(x)| < \sum_{1 \leq i \leq k} l(i)$ (all elements of\textit{ E(x)} are natural numbers representable using $k$ symbols or less which means the number of possible terms using $k$ symbols or less has to be strictly larger than order of $E(x)$ (some terms of length $k$ don't evaluate to an actual integer, ex: \textit{‘***....***’})).//
Therefore, $|E(x)| < \sum_{1 \leq i \leq k} l(i) < \sum_{1 \leq i \leq k} a^{i} = \mathcal{O}(x^{1-\epsilon})$\\
\textbf{Note}: This theorem is a generalization of \textbf{Theorem 1.4} obtained by Gnang, Radziwill and Sanna in [3].
\\
\\
\\
\\
\textit{\textbf{Theorem 1.5:} $\forall k,k'$\textbf{ s.t.} $k<k'$, \textbf{with corresponding ugly numbers} $n_u,n'_u$\textbf{ under an operational set s.t.} $S\in O$, $n_u<n'_u$\textbf{.}}\\
\textbf{Proof:} Let $n'_u$ be the smallest ugly number s.t. $c_O(n'_u) = k'$ and $\exists$ ugly number $n_u$ with $c_O(n_u) = k$ satisfying $k<k'$, $n'_u<n_u$. \textbf{Question:} What is $c_O(n'_u - 1)$ w.r.t $k$.\\
If $c_O(n'_u - 1) < k < k'$, then I can use a successor function on the term evaluating to $n'_u - 1$ to get $n'_u$, allowing me to write $c_O(n'_u) \leq k < k'$ which is a contradiction by our construction.
\\
If $c_O(n'_u - 1) = k'' > k$, $\exists n''\leq n'_u - 1 < n'_u < n$, s.t. $c_O(n'') = k''>k \implies \exists n''_u\leq n'_u - 1 < n'_u < n_u$, s.t. $c_O(n''_u) = k'' > k$. This is a contradiction by construction as $n'_u$ was supposed to be the smallest ugly number fulfilling that property.
\\
If $c_O(n'_u - 1) = k$, $\exists n'' \leq n'_u - 1 < n'_u < n_u$, s.t. $c_O(n'') = k \implies \exists n''_u \leq n'_u - 1 < n'_u < n_u$, s.t. $c_O(n''_u) = k$. But $k = c_O(n_u)$. Therefore $c_O(n'_u - 1)$, a quantity we know has to be finite, is neither smaller, equal or larger than $k$, another quantity we know to be finite. This is the final contradiction. Hence, our initial hypothesis was wrong and there exists no $n'_u$ that can satisfy our requirements.
\\
\textbf{Corollary: $\forall$ consecutive natural numbers $k_1 < k_2 < .. < k_m < ..$, the corresponding ugly numbers are in order $n_1 < n_2 < .. < n_u < ..$ Therefore, $\forall n< n_{u}$, $c_O(n_u) = k+1$, $c_O(n)\leq k$.}
\\
\\
\textbf{\textit{Theorem 1.6: $\forall O$, s.t. $S\in O$, if $n + 1$ is ugly and $c_{O}(n+1) = k+1$, then $n+1$ has a minimal representation of the form $S'r'$ where $v(r) = n$}}\\
\textbf{\textit{Proof:}} $c\textsubscript{\textit{O}}(n) \leq k$ by \textbf{Corollary} to \textbf{\textit{Theorem 1.5}}. If $c\textsubscript{\textit{O}}(n) < k$, then $c\textsubscript{\textit{O}}(n+1) \leq k$ as $v(S'r') = n+1$ and $S'r'$ uses at most $(k - 1) + 1$ symbols. This is a contradiction by construction. Therefore minimal representation has to have length $k+1 \implies S'r'$ is a minimal representation of our ugly number.

\pagebreak
\section{Results on $\{1, S, *\}$ and $\{1, S, +, *\}$}
\textit{\textbf{Theorem 2.1: Let} $O = \{1, S, *\}$\textbf{.} $c\textsubscript{O}(n) \geq (\gamma + 1) \log(n)/\log(\gamma) - 1$ \textbf{where}\\
$\gamma$ = $\mathbf{argmax\textsubscript{w}} (\log(w)/(w + 1)) ≈ 3.59$.}
\\
\textit{\textbf{Proof:}} Fix a value $k$ for the length. Let $r$ be the term of length $k$ with maximal value.
\\
It is never worthwhile to have a $*$ inside the scope of an $S$, because for any subterms of value x, y $>$ 1, $v(S*xy)<v(*Sxy)$. Therefore, since multiplication is associative, we can assume that $r$ consists of $q − 1$ occurrence of $∗$ followed by $q$ terms $x_{1} . . . x_{q}$ each of which has of the form $'SS . . . S1'$, for varying number of $S$.
\\
Therefore we have $v(r) = \prod x\textsubscript{q}$. We need to maximize this term keeping the constraint that the total complexity value of the term remain equal to $k$. Therefore, to make the problem easier, we solve the problem in the reals rather than in the integers, thus giving us a weaker lower bound for complexity of $n$.
\\
We use the classical result that if the sum of a number of real-valued terms is fixed, the product is maximized when the terms are all equal. Putting $x_{1} = x_2 = x_3 ... = x_{q} = w \implies v(r) = w^{q}$ and $qw + q - 1= k$; thus $q = (k + 1)/(w + 1)$.  We need to solve:
\\
$\mathbf{argmax}_{w} (v(r))=\mathbf{argmax}_{w} (w^{(k+1)/(w+1)}) = \mathbf{argmax}_{w} (log(w)/(w+1) )= \gamma$ \\($k+1$ is constant and hence can be ignored after we take log on both sides and consider $\mathbf{argmax}_{w}$)
\\
Therefore $n = v(r) \leq \gamma^{(k+1)/(\gamma+1)}$. Taking the log on both sides and re-arranging gives us the quoted result.
\\
\textbf{Corollary: Let $ O = \{1, S, +, *\}$. $c_{O}(n) \geq (\gamma + 1) \log(n)/\log(\gamma) - 1$ where\\
$\gamma = \mathbf{argmax}_{w} (\log(w)/(w + 1)) ≈ 3.59$. \\
(Like S, + is not used for the maximal number)}
\\
\\
\textit{\textbf{Theorem 2.1 (Stronger statement): Let $O = \{1, S, *\}$. $c_{O}(n) \geq 5\log(n)/\log(4) - 1$}}
\\
\textit{\textbf{Proof:}} Fix a value $k$ for the length. Let $'r'$ be the term of length $k$ with maximal value. $'r'$ has general representation of form $*'a'*'b'*'b'....(c $ copies$)....*'b''b'\equiv a.b^c$ s.t. $a + bc + c = k$. Let $*'b'*'b'...*'b''b'$ use $k_1$ symbols. $f(w) = log(w)/(w + 1)$ is a monotonic decreasing function beyond $\gamma$. Further, $f(1), f(2), f(3) < f(4)$. Therefore, for $m$ s.t. $5m - 1 = k_1, 4^m \geq b^c$ (following the constraint argument developed earlier). The same argument yields, for $m_a$ s.t. $5m_a - 1 = a, 4^{m_a} \geq a$. Now, for $m_k$ s.t. $m_k = m_a + m$, $4^{m_k} = 4^{m_a + m} \geq a.b^c = v(r) \implies m_k \geq \log_4(r)$. Combining these findings with \textbf{Theorem 1.3} gives us our result.\\
\textbf{Corollary: Under $O = \{1, S, *\}$, numbers of the form $4^k, k\in \textbf{N}$ are the only numbers to achieve the lower bound.}\\
\textbf{Note:} An equivalent formulation of the corollary was shown in [4].
\\
\\
\\
\textit{\textbf{Theorem 2.2: Let $O = \{1, S, *\}$. $c_{O}(n) \leq 8\log(n)/\log(4) + 2$}}
\\
\textit{\textbf{Proof:}} $\forall n$ s.t. $4^{k} \leq n \leq 4^{k+1}$ for some $k \geq 1$, $\exists a_{i}$ satisfying $1 \leq a_{i} \leq 3$, s.t. $n = \sum_{1 \leq i \leq k} a_i4^i$
\\Re-arranging, $n = a_{0} + 4(a_i + 4(a_2 + 4(...(a_{k - 1} + 4a_{k}))))$. Using $S$ to express the $a_{i}$, and counting the amount of symbols used in the most inefficient scenario (if $a_{i} = 3, \forall i$, we could factor $3$ out. For $i = k$, if we let $a_{i} = 2$, we count for the ensemble using the most number of symbols), 
\\
$c_{O}(n) \leq 8k + 2$. Using $4^{k} \leq n \leq 4^{k+1}$, $c_{O}(n) \leq 8\log(n)/\log(4) + 2$
\\
\\
\textbf{\textit{Proposition 2.1: Let $O = \{1, S, *\}$ or $\{1, S, +, *\}$. $'M_O(k)'$ can never contain any sub-terms evaluating to $6, 7$ and more than 1 sub-term evaluating to $2$ or $5$. There can be at most $4$ sub-terms evaluating to $3$.}}\\
\textbf{\textit{Proof:}} If a factor in a term can be replaced by a higher factor in the term, the term is not a maximal representation. As already discussed, $S$ and $+$ outside scope of $*$ don't generate maximal values. $7 \equiv$ $'SSSSSS1' \implies$ we can replace it by $'*SS1SS1'$. This holds true for any set that contains $1, S, *$ and has $c_{O}(7) = 7$. Let X be any expression that computes to a value.\\
There can't be more than $1$ factor of $6$ as $'X*SSSSS1SSSSS1'$ is better written as $'X*SS1*SSS1SSS1'$. Therefore at most $1$ factor of $6$ survives. If $X$ contains even $1$ factor of $5$, we can write $'*SSSSS1SSSS1'$ as $'*SSS1*SS1SS1'$. If $X$ contains no $5$, it contains factors of $2,3$ or $4$. $'*6A'$ where $A$ is a stand in for any of $2, 3$ or $4$ is better written as $'*5(A + 1)'$. Therefore no factor of $6$ survives in a maximal representation.\\
There cannot be more than one $2$ because we can always pair $'*S1S1'$ as $'SSS1'$ and use the one remaining symbol to increase the value of our term. Further, if there is $4$ or $5$, we can replace $2$ by $3$ and decrease the value of either of them by $1$. Further if there are two factors of $3$, we can rewrite $'X*S1*SS1SS1'$ better as $'X*SSS1SSSS1'$. As $'*S1SS1' \equiv$ $'SSSSS1'$, the only maximal representation containing $2$ is $2 \equiv$ $'S1'$.\\
If there are more than one factors of $5$,  $'X*SSSS1SSSS1'$ can be better written as $'X*SS1*SS1SS1'$. Further, if there are at least two * involved, $5$ is not involved in the maximal representation. This is because $'*SSSS1SS1'$ is better written as $'*SSS1SSS1'$. Therefore, if even one factor of $3$ is present, a factor of $5$ cannot be present. If it is all $4$s, $'*SSS1*SSS1SSSS1'$ is better written as $'*SS1*SS1*SS1SS1'$. Therefore only maximal representations containing $5$ are $5 \equiv$ $'SSSS1'$ and $20 \equiv$ $'*SSS1SSSS1'$. \\
$v(*SS1*SS1*SS1*SS1SS1) < v(*SSS1*SSS1*SSS1SSS1) \implies$ there can't be more than four factors of 3.\\
\textbf{Note:} The procedures described above for 6 and 7 generalize to sets composed of $1,S,*$ and other arbitrary order hyper-operations. A case by case approach is required for 2, 3, 4, 5.\\
\textbf{Note:} A discussion of the results on $'S1'$ and $'SSSSS1'$ was done in [4].\\
\textbf{Theorem 2.3: Let $k = 5m - 1 - r \geq 11$ where $m$ is the number such that 4 $\geq r\in \mathbf{N} \geq0$. Then $v(M_O(k)) = 3^{r}4^{m - r}$}.\\
\textbf{Proof:} $\forall k\geq$ 11, only factors of 3 and 4 are allowed and there can never be more than 4 factors of $3$.\\
\\
\section{Observations}
As a way of testing results and gaining insight, the computing resources at the Courant Institute were used to generate the minimal representations for the first $4.5*10^6$ numbers. This section presents some relevant observations and patterns from the data, with minor comments wherever possible. \\
\\
\textbf{Observation 3.1: Let $O = \{1, S, +, *\}$. $c_{O}(n) < (c_O(n_u) + 1)(a + 1) - 2$ where $a\in N$ is s.t. $n \leq a^{a}$ and $n_u$ is the largest ugly number satisfying $c_O(n_u) \leq \lceil a \rceil$} (plotted in \textbf{red}).
\\
\\
\textbf{Observation 3.2: Let $O = \{1, S, +, *\}$. $\forall n \leq 4500000, c_{O}(n) \leq\left \lceil{5log(n)/log(4) + a + 1}\right \rceil$ where $a$ is as defined above.} \textbf{The function is plotted in green below.}
\\
\\
\textbf{Observation 3.3: Let  $O' = \{1, S, +, *\}$ and  $O = \{1, S, *\}$.} $\forall n \leq 4500000, c\textsubscript{O'}(n) = c\textsubscript{O}(n).$
The graph below plots the complexity value of the first 4.5 million natural numbers, along with the upper and lower bounds we found for those values. The most surprising fact borne out of our computational efforts was that $O'$ and $O$ had no difference in minimal representations for the first 4.5 million numbers.\\
\\
\textbf{Observation 3.4: If $k > k'$, then $k$ is a minimal representation more often than $k'$.}
From \textbf{Defn} 1 to 4, it is clear that $k$ has more representations than $k'$ (we can make all the possible representations made by m with n too, with the additional units serving as fixed placeholders at the start or end. Then we can permute those placeholders to generate more representations). It seems therefore intuitive that $k$ also has more \textbf{minimal} representations to its name than $k'$. We expect and observe monotonic increasing behavior in the number of natural numbers having the same length of minimal representation as number of symbols allowed is increased.
\\
\begin{figure}[H]
 \begin{center}
  \includegraphics[width=18cm]{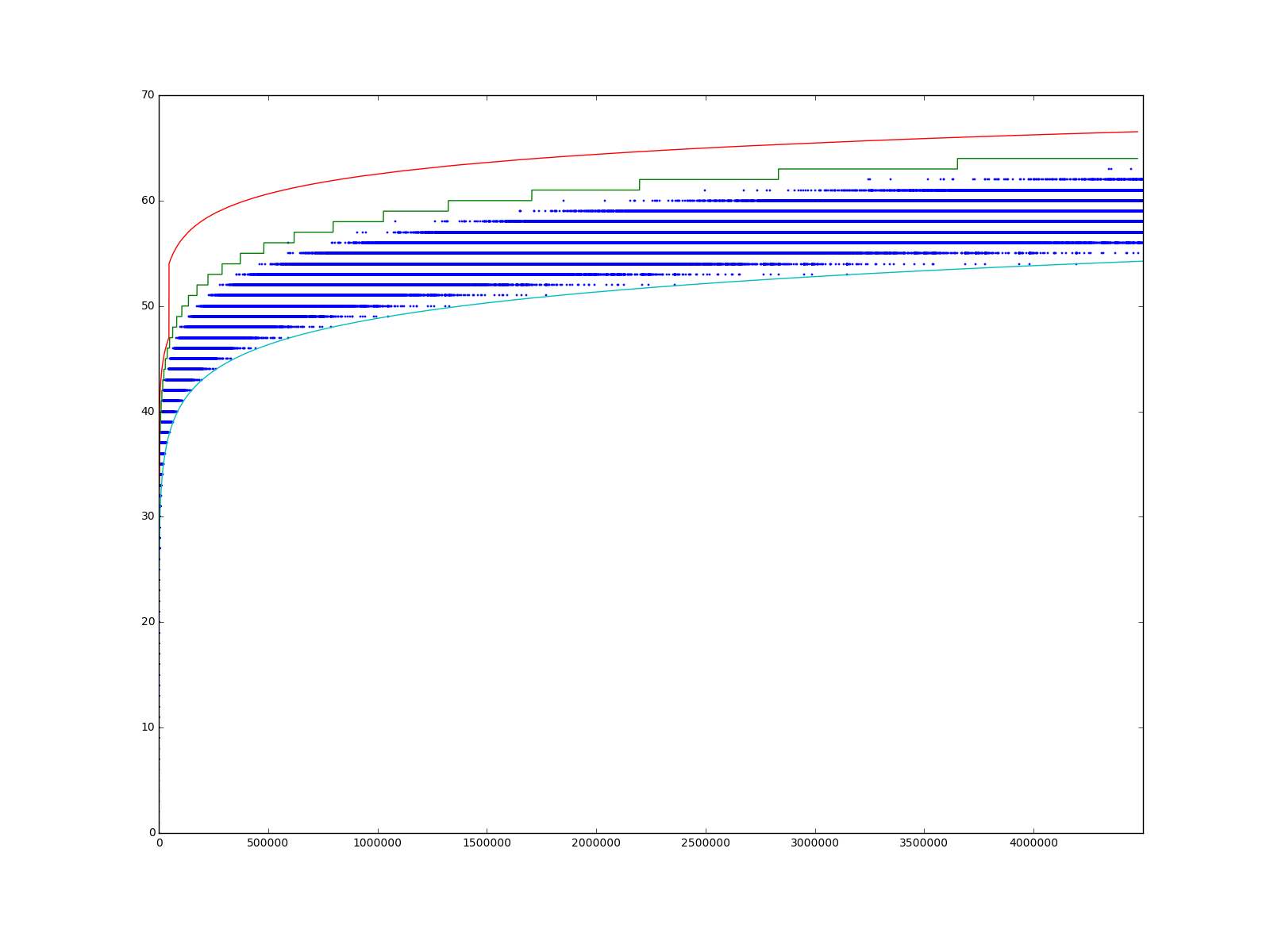}
  \caption{Bounds on complexity of Minimal Representations under $O, O'$}
  \label{fig:1S+}
 \end{center}
\end{figure}

\begin{figure}[H]
 \begin{center}
  \includegraphics[width=18cm]{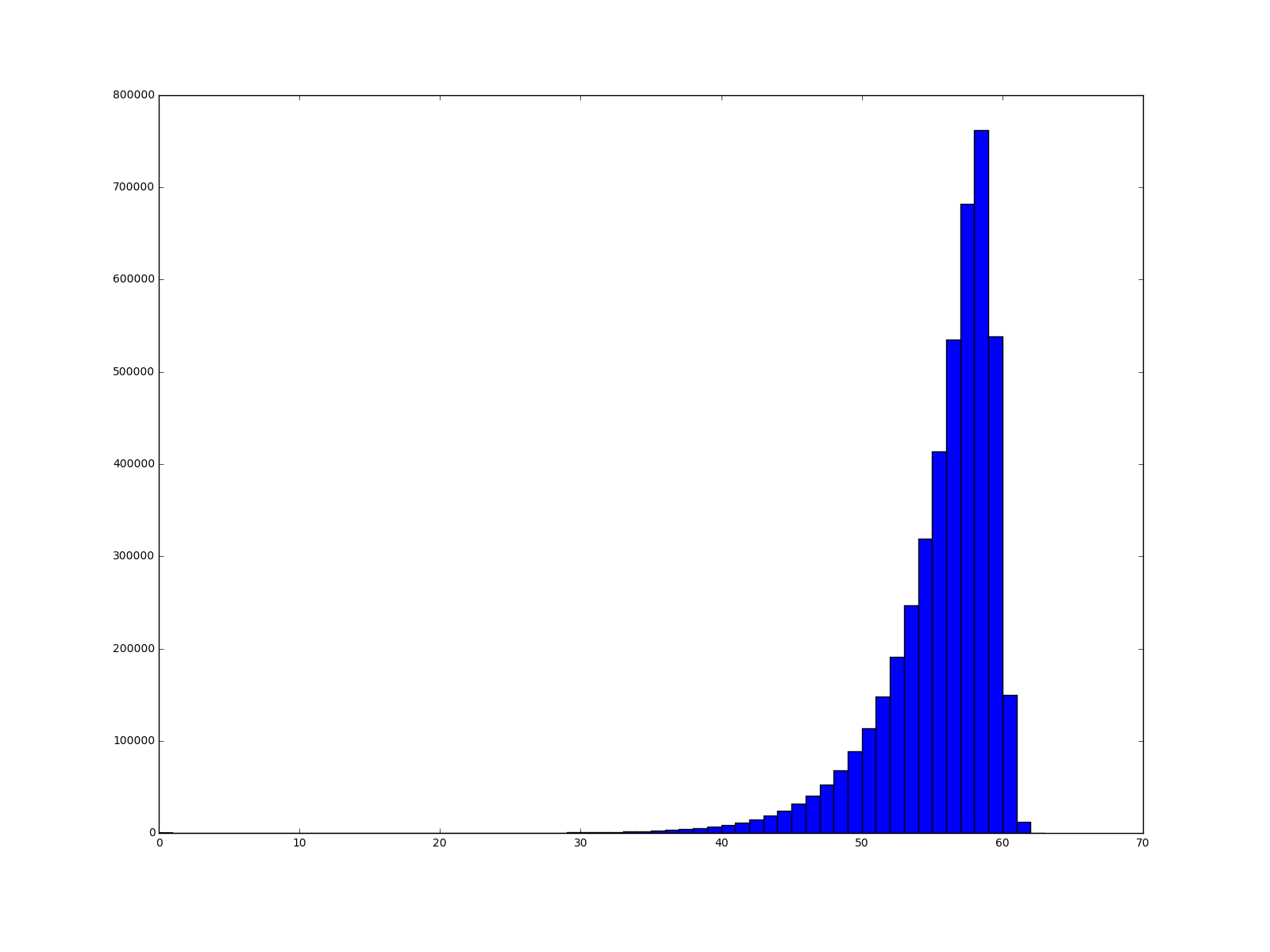}
  \caption{Numbers with complexity $k$ v.s. $k$ under \{1, S, +, *\}}
  \label{fig:Freqanal1Smulti}
 \end{center}
\end{figure}

\begin{table}[]
\centering
\caption{Ugly Numbers for $7< k\leq 63$ under $O = \{1, S, *\}$}
\label{First}
\begin{tabular}{llll}
$n_u$ & Minimal Representation & $c_O(n_u)$ & Primality \\
& & & \\
10 & S*SS1SS1 & 8 & Not Prime \\
11 & SS*SS1SS1 & 9 & Prime \\
14 & SS*SS1SSS1 & 10 & Not Prime \\
19 & S*S1*SS1SS1 & 11 & Prime \\
22 & SS*SSS1SSSS1 & 12 & Not Prime \\
23 & SSS*SSS1SSSS1 & 13 & Prime \\
38 & SS*SS1*SS1SSS1 & 14 & Not Prime \\
43 & S*S1S*SSS1SSSS1 & 15 & Prime \\
58 & S*SS1S*S1*SS1SS1 & 16 & Not Prime \\
59 & SS*SS1S*S1*SS1SS1 & 17 & Prime \\
89 & S*S1*SSS1SS*SS1SS1 & 18 & Prime \\
107 & SS*SS1*SSSS1SSSSSS1 & 19 & Prime \\
134 & SS*SS1*SSS1SS*SS1SS1 & 20 & Not Prime \\
167 & SS*SS1S*S1*SS1*SS1SS1 & 21 & Prime \\
179 & SSS*SSS1*SSS1SS*SS1SS1 & 22 & Prime \\
263 & SSS*SSS1S*SSS1*SSS1SSS1 & 23 & Prime \\
347 & SS*SSSS1S*SSS1S*SSS1SSS1 & 24 & Prime \\
383 & SSS*SSS1*SSSS1S*S1*SS1SS1 & 25 & Prime \\
537 & S*S1*SSS1SSS*SSS1*SSS1SSS1 & 26 & Not Prime \\
713 & SS*SS1*SS1S*S1*SS1S*SS1SSS1 & 27 & Not Prime \\
719 & SS*SS1S*S1*SSSSSS1S*SSS1SSS1 & 28 & Prime \\
1103 & SSS*SSS1*SSSS1S*S1*SS1*SS1SS1 & 29 & Prime \\
1319 & SSS*SSS1S*SSS1S*SS1*SS1*SS1SS1 & 30 & Prime \\
1439 & SS*SS1SS*SS1*SS1S*SSS1S*SS1SSS1 & 31 & Prime \\
2099 & SSS*SSS1*SSS1S*S1S*SSS1*SSS1SSS1 & 32 & Prime \\
2879 & SS*SS1*SSSSSS1SS*SS1*SS1*SS1SSSS1 & 33 & Prime \\
3833 & SSS*SSSS1S*SS1*SS1*SSSS1S*SSS1SSS1 & 34 & Prime \\
4283 & SSS*SSSS1S*SS1*SS1*SSSS1S*S1*SS1SS1 & 35 & Prime \\
5939 & SSS*SSS1*SSS1S*S1*SSSS1S*SS1*SS1SSS1 & 36 & Prime \\
6299 & SSS*SSS1S*S*SS1SSS1S*S1*SS1*SSS1SSSS1 & 37 & Prime \\
9059 & SSS*SSS1*SSS1S*SSSS1S*SSS1S*SS1*SS1SS1 & 38 & Prime \\
12239 & SSS*SSS1*SSSSSS1S*SSS1S*SS1*SS1*SS1SSS1 & 39 & Prime \\
15118 & SS*SSS1S*S1S*SSS1*SSS1S*SS1*SS1S*SS1SSS1 & 40 & Not Prime \\
19079 & SSSS*SSSS1*SSSS1*SSSSSS1S*SS1*SS1*SS1SSS1 & 41 & Prime \\
23039 & SSSS*SSSS1*S*SSS1SSS1S*S1*SS1*SS1*SS1SSSS1 & 42 & Prime \\
26459 & SSSS*SSSS1*SS*SS1SS1S*S1*SS1*SSS1*SSS1SSSS1 & 43 & Prime \\
44879 & SS*SS1S*S1*SS1*SS1*SS1S*S1*SS1S*SS1*SS1SSSS1 & 44 & Prime \\
49559 & SSSS*SSSS1*SS*SS1SS1S*SS1*SS1*SSS1*SSSS1SSSS1 & 45 & Prime \\
66239 & SSS*SSS1S*S1*S*SSS1SSS1S*S1*SS1*SS1*SS1*SS1SS1 & 46 & Prime \\
78839 & SSS*SSS1S*SSS1*S*SS1SSS1S*S1*SS1*SS1S*SSS1SSSS1 & 47 & Prime \\
98999 & SS*SS1SS*S*SSS1SSS1S*SSS1*SSSS1S*S1*SS1*SSS1SSS1 & 48 & Prime \\
137339 & SSSS*SSSS1*SS*SS1SS1S*SS1*SSS1*SSS1*SSS1S*SS1SSS1 & 49 & Prime \\
172583 & SSS*SSS1S*S1*SSS1S*SSS1*SSS1S*SS1*SSS1S*SS1*SS1SS1 & 50 & Prime \\
228479 & SSSS*SSSS1*SSSS1*S*SS1SSS1S*S1*SS1*SS1*SS1S*SS1SSS1 & 51 & Prime \\
280223 & SSS*SSS1*SSSS1S*S1*SSSS1S*S1*SSS1*SSSS1*SSSS1SSSSSS1 & 52 & Prime \\
355679 & SSS*SSS1S*S*SS1*SS1SSSS1S*SS1*SSS1S*S1*SSS1*SSS1SSSS1 & 53 & Prime \\
460079 & SSSS*SSSS1*SSSS1*SSSSSS1S*SS1*SS1*SSS1S*S1*SS1*SS1SSS1 & 54 & Prime \\
590398 & SSS*SSSS1*S*SS1SSS1SS*SS1*SS1S*SS1*SS1*SSS1S*SS1*SS1SS1 & 55 & Not Prime \\
590399 & SSSS*SSSS1*S*SS1SSS1SS*SS1*SS1S*SS1*SS1*SSS1S*SS1*SS1SS1 & 56 & Prime \\
907199 & SSS*SSS1S*S1S*S*SS1*SS1S*SS1SSS1S*SS1*SSS1*SSS1*SSS1SSSS1 & 57 & Prime \\
1081079 & SSSS*SSSS1*SSSS1S*S1S*SSS1*SSSS1S*S1*SS1*SS1*SS1*SSS1SSSS1 & 58 & Prime \\
1650983 & SSS*SSS1*SSSS1S*SS1*SS1*SSS1S*SS1S*SSSSSS1S*SS1*SS1*SS1SSS1 & 59 & Prime \\
1851119 & SSSS*SSSS1S*S*SS1*SS1SSS1S*SSSS1S*SSS1*SSS1*SSSS1*SSSS1SSSS1 & 60 & Prime \\
2497499 & SSSS*SSSS1*SSSSSS1*S*SS1SSS1S*SSS1*SSS1*SSSSSS1S*SS1*SSS1SSS1 & 61 & Prime \\
3243239 & SSSS*SSSS1*SS*SS1*SS1SSSS1S*SS1*SSS1*SSSS1*SSSS1S*SS1*SS1SSSS1 & 62 & Prime \\
4344479 & SSSS*SSSS1*SSSS1S*SS1S*SSSS1S*SSS1*SSS1*SSS1S*SS1*SS1*SSS1SSSS1 & 63 & Prime \\
\end{tabular}
\end{table}

\begin{table}[]
\centering
\caption{$M_O(k)$ for $1\leq k \leq 54$ under $O = \{1, S, *\}$}
\label{Last}
\begin{tabular}{ll}
$k$ & $M_O(k)$\\
\\
1 & 1 \\
2 & S1 \\
3 & SS1 \\
4 & SSS1 \\
5 & SSSS1 \\
6 & SSSSS1 \\
9 & *SS1SS1 \\
12 & *SS1SSS1 \\
16 & *SSS1SSS1 \\
20 & *SSS1SSSS1 \\
27 & *SS1*SS1SS1 \\
36 & *SS1*SS1SSS1 \\
48 & *SS1*SSS1SSS1 \\
64 & *SSS1*SSS1SSS1 \\
81 & *SS1*SS1*SS1SS1 \\
108 & *SS1*SS1*SS1SSS1 \\
144 & *SS1*SS1*SSS1SSS1 \\
192 & *SS1*SSS1*SSS1SSS1 \\
256 & *SSS1*SSS1*SSS1SSS1 \\
324 & *SS1*SS1*SS1*SS1SSS1 \\
432 & *SS1*SS1*SS1*SSS1SSS1 \\
576 & *SS1*SS1*SSS1*SSS1SSS1 \\
768 & *SS1*SSS1*SSS1*SSS1SSS1 \\
1024 & *SSS1*SSS1*SSS1*SSS1SSS1 \\
1296 & *SS1*SS1*SS1*SS1*SSS1SSS1 \\
1728 & *SS1*SS1*SS1*SSS1*SSS1SSS1 \\
2304 & *SS1*SS1*SSS1*SSS1*SSS1SSS1 \\
3072 & *SS1*SSS1*SSS1*SSS1*SSS1SSS1 \\
4096 & *SSS1*SSS1*SSS1*SSS1*SSS1SSS1 \\
5184 & *SS1*SS1*SS1*SS1*SSS1*SSS1SSS1 \\
6912 & *SS1*SS1*SS1*SSS1*SSS1*SSS1SSS1 \\
9216 & *SS1*SS1*SSS1*SSS1*SSS1*SSS1SSS1 \\
12288 & *SS1*SSS1*SSS1*SSS1*SSS1*SSS1SSS1 \\
16384 & *SSS1*SSS1*SSS1*SSS1*SSS1*SSS1SSS1 \\
20736 & *SS1*SS1*SS1*SS1*SSS1*SSS1*SSS1SSS1 \\
27648 & *SS1*SS1*SS1*SSS1*SSS1*SSS1*SSS1SSS1 \\
36864 & *SS1*SS1*SSS1*SSS1*SSS1*SSS1*SSS1SSS1 \\
49152 & *SS1*SSS1*SSS1*SSS1*SSS1*SSS1*SSS1SSS1 \\
65536 & *SSS1*SSS1*SSS1*SSS1*SSS1*SSS1*SSS1SSS1 \\
82944 & *SS1*SS1*SS1*SS1*SSS1*SSS1*SSS1*SSS1SSS1 \\
110592 & *SS1*SS1*SS1*SSS1*SSS1*SSS1*SSS1*SSS1SSS1 \\
147456 & *SS1*SS1*SSS1*SSS1*SSS1*SSS1*SSS1*SSS1SSS1 \\
196608 & *SS1*SSS1*SSS1*SSS1*SSS1*SSS1*SSS1*SSS1SSS1 \\
262144 & *SSS1*SSS1*SSS1*SSS1*SSS1*SSS1*SSS1*SSS1SSS1 \\
331776 & *SS1*SS1*SS1*SS1*SSS1*SSS1*SSS1*SSS1*SSS1SSS1 \\
442368 & *SS1*SS1*SS1*SSS1*SSS1*SSS1*SSS1*SSS1*SSS1SSS1 \\
589824 & *SS1*SS1*SSS1*SSS1*SSS1*SSS1*SSS1*SSS1*SSS1SSS1 \\
786432 & *SS1*SSS1*SSS1*SSS1*SSS1*SSS1*SSS1*SSS1*SSS1SSS1 \\
1048576 & *SSS1*SSS1*SSS1*SSS1*SSS1*SSS1*SSS1*SSS1*SSS1SSS1 \\
1327104 & *SS1*SS1*SS1*SS1*SSS1*SSS1*SSS1*SSS1*SSS1*SSS1SSS1 \\
1769472 & *SS1*SS1*SS1*SSS1*SSS1*SSS1*SSS1*SSS1*SSS1*SSS1SSS1 \\
2359296 & *SS1*SS1*SSS1*SSS1*SSS1*SSS1*SSS1*SSS1*SSS1*SSS1SSS1 \\
3145728 & *SS1*SSS1*SSS1*SSS1*SSS1*SSS1*SSS1*SSS1*SSS1*SSS1SSS1 \\
4194304 & *SSS1*SSS1*SSS1*SSS1*SSS1*SSS1*SSS1*SSS1*SSS1*SSS1SSS1 \\
\end{tabular}
\end{table}

\pagebreak

\pagebreak
\section{Tackling higher order hyper-operations}
\textbf{{Theorem 4.1: Let $O = \{1, S, \wedge \}$. $c_{O}(n) \geq 4\log*{3}(n) - 1$}}\\
Proof: Using the same arguments as in \textbf{Thm 2.1}, we find that $a^{a^{a^{a^{.^{.^{(b\ times)}}}}}} \equiv a\uparrow\uparrow{b}$ is the most efficient way of building up number. $n = a\uparrow\uparrow{b} \implies b = \log*\textsubscript{a}n \implies$ function to minimize is $(x+1)(1+\log(1/\log x)/\log x) \implies a = 3 \implies c_{O}(n) \geq 4\log*_{3}(n) - 1.$ 
\\

\section{Acknowledgements}
I would like to express my gratitude for Professor Ernest Davis of Courant Institute, NYU, for stimulating discussions on the mathematical and computational objectives of the project, along with his insight into the problems and the methods to tackle these questions. I am also grateful for the computational inputs by Mr. Arvi Gjoka of Courant Institute, primarily in compiling the minimal representations of the first 4.5 million numbers. A similar thanks to Mr. Adolfo Holguin of Courant Institute, for his help on proof-reading this work. I would also like to thank the Dean's Undergraduate Research Fund (DURF) committee at NYU, for deciding to fund this research. Special thanks to Mr. Carlo Sanna of Università degli Studi di Torino for a fruitful discussion on the problem and to the $20^{th}$ International Workshop for Young Mathematicians, for hosting us both for a lecture.\\
\\
Lastly, and most importantly, I would be remiss if I did not mention Randall Kayser, Iraj Eshghi of New York University, Daniel Polin of UC Davis, Sreela Das of McGill University and George Wong of University of Illinois at Urbana–Champaign. The inspiration for this project came from a game we played while we were all studying in New York University, of making numbers by using arithmetical operations on $1$ and quite naturally wanted to know the fastest ways to do it. In particular, Mr. Polin was the first of us to obtain a lower bound on the complexity when $*$ is the highest allowed hyperoperation. The stronger statement of \textbf{\textit{Theorem 2.1}}, although obtained from different considerations and in a different setting, is equivalent to his own findings (and to the findings of Mr. Matos, which are given in [4]). I would also like to thank Mr. Wong for our discussions on the best way to represent arithmetical computations and general discussions on how we should approach the project. A special thanks is also deserved for Mr. Matos for his work in [4].
\section*{References}
\begin{enumerate}
\item R. L. Goodstein (Dec 1947). "Transfinite Ordinals in Recursive Number Theory". Journal of Symbolic Logic. 12 (4): 123–129
\item Kolmogorov, Andrey (1963). "On Tables of Random Numbers". Sankhyā Ser. A. 25: 369–375. MR 0178484.
\item Edinah K. Gnang, Maksym Radziwiłł, Carlo Sanna (July 2015). "Counting Arithmetical Formulas". European Journal of Combinatorics: Vol 47, Pages 40-53
\item Armando P. Matos (July 2015). "Kolmogorov complexity in multiplicative arithmetic". DCC-Faculdade de Ciencias da Universidade do Porto, Technical Report
\end{enumerate}

\end{document}